\documentclass[a4paper,10pt]{article}
\usepackage[utf8]{inputenc}
\usepackage{amsthm}
\usepackage{amsmath,amssymb,amsfonts}
\usepackage[dvipsnames,usenames]{color}
\usepackage[bitstream-charter]{mathdesign}
\usepackage[utf8]{inputenc}
\usepackage{tikz}
\usepackage{caption}
\usepackage{subcaption}
\usepackage{enumerate}
\usepackage{geometry}
\usepackage{authblk}
\usepackage{tikz-cd}
\title{Planar cubic graphs of small diameter}
\author[1,2]{Kolja Knauer}
\author[3]{Piotr Micek}

\affil[1]{{\small Aix-Marseille Univ, Universit\'e de Toulon, CNRS, LIS, Marseille, France}}\affil[2]{{\small Departament de Matem\`atiques i Inform\`atica, Universitat de Barcelona (UB), Barcelona, Spain}}

\affil[3]{{\small Jagiellonian University,
Faculty of Mathematics and Computer Science, Theoretical Computer Science Department, Poland}}

\date{}

\begin{document}

\maketitle

\begin{abstract}
Cubic planar $n$-vertex graphs with faces of length at most $6$, e.g.~fullerene graphs, have diameter in $\Omega(\sqrt{n})$. It has been suspected, that a similar result can be shown for cubic planar graphs with faces of bounded length. This note provides a family of cubic planar $n$-vertex graphs with faces of length at most $7$ and diameter in ${O}(\log n)$, 
thus refuting the above suspicion.
\end{abstract}

%\section{The construction}
A \emph{fullerene graph} is a cubic $3$-connected planar graph all of whose faces are of length $5$ or $6$. They model fullerene molecules, i.e., molecules consisting only of carbon atoms other than graphite and diamonds. 
The search for correlations of graph invariants and chemical stability motivates the study of fullerenes in chemical graph theory, see e.g.~\cite{AKS16,SWA15}. 

The diameter of a planar cubic graph can be logarithmic in the number of vertices. However, the diameter of a fullerene graph on $n$ vertices is at least $\frac{1}{6}\sqrt{24n-15}-\frac{1}{2}$, as proved in~\cite{ADKLS12}. Indeed, the proof of this lower bound extends to planar cubic graphs whose faces are of size at most~$6$. It has been suspected, that a similar lower bound holds for the diameter of cubic planar graphs with faces of bounded length, see~\cite{ADKLS12,AKS16}.
This is false as shown by the following family $\{G_k\mid k\geq 2\}$.

\begin{figure}[ht]
\centering
\includegraphics[width=.6\textwidth]{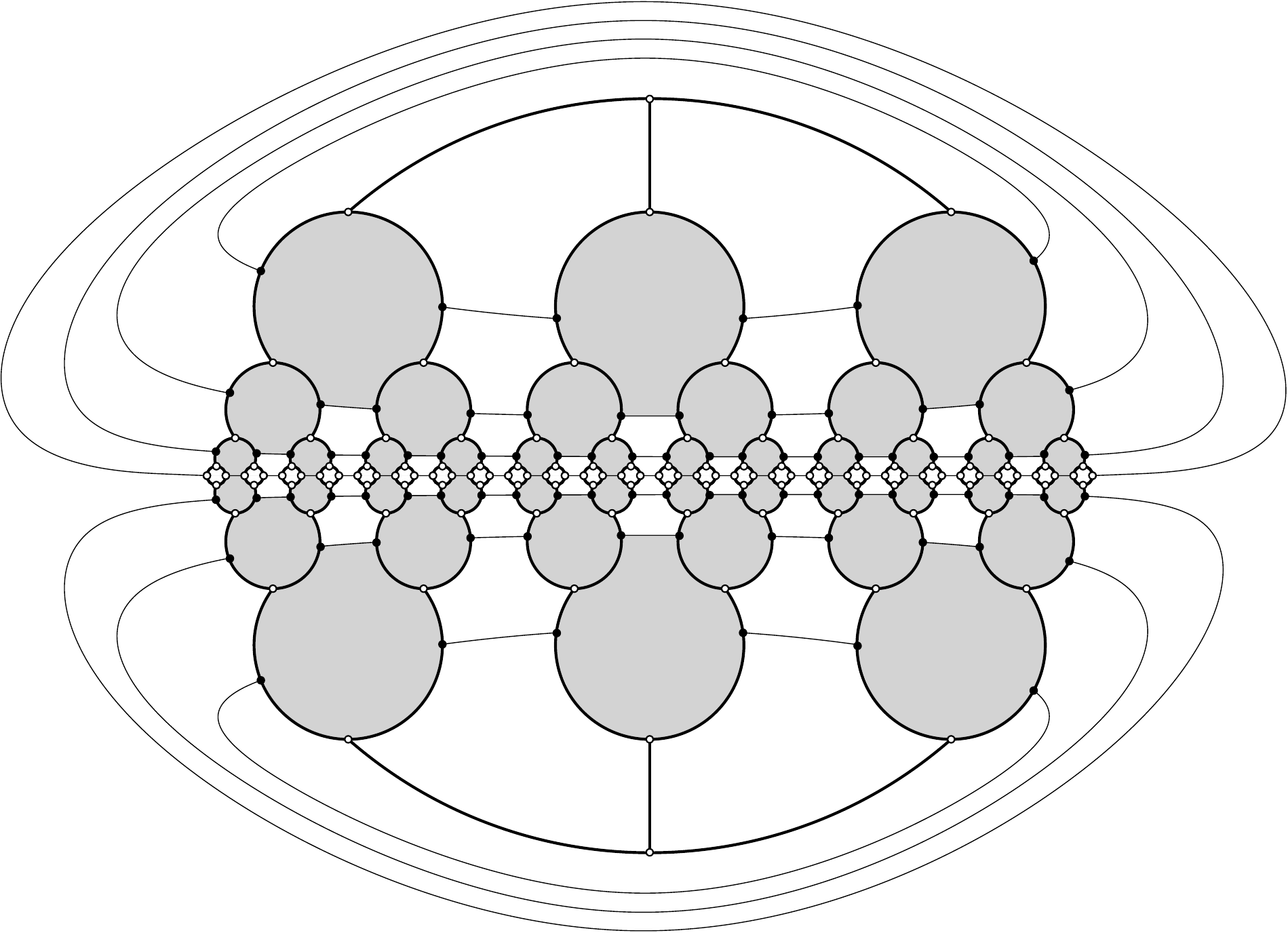}
\caption{The graph $G_5$. The vertices of the ternary trees are white and the original tree edges are bold. The gray faces are of length $7$, all the others are shorter.}
\label{fig:fullgraph}
\end{figure}

Let $k\geq 2$. Denote by $T_k$ the rooted tree
all of whose internal vertices are of degree $3$ and all leaves are at distance $k$ from the root.
To define $G_k$, first take two copies of $T_k$ and glue them by identifying corresponding leaves, see Figure~\ref{fig:fullgraph}. Now, subdivide all edges except those incident to leaves or roots. In Figure~\ref{fig:fullgraph}, the new vertices are black, while the original vertices of the trees are white. 
Fix a plane embedding of the obtained graph such that both roots are on the outer face. 
Note that, for each tree the vertices in the tree of a given distance $d$ to its root are naturally ordered from left to right.
Note that all the subdivision vertices and all the leaves are exactly the vertices at even distance to the root of their tree.
Now, for each of the two trees and for every even integer $d$ with $2\leq d \leq 2k-2$, let $v_0,\ldots,v_{\ell}$ be the vertices at distance $d$ from the root in the tree ordered from left to right.
Note that $\ell$ is odd as $d\geq2$.
Put edge $v_iv_{i+1}$ for every even $i$ 
(where $i+1$ is taken modulo $\ell+1$).
This completes the definition of $G_k$.

The added edges between leaves and subdivision vertices did not destroy planarity, 
see Figure~\ref{fig:fullgraph}.
After gluing the two trees together and subdividing the edges the only vertices of degree $2$ were the leaves and the subdivision vertices.
Each of them got one new incident edge, thus $G_k$ is cubic.

\begin{figure}[ht]
\centering
\includegraphics[width=.65\textwidth]{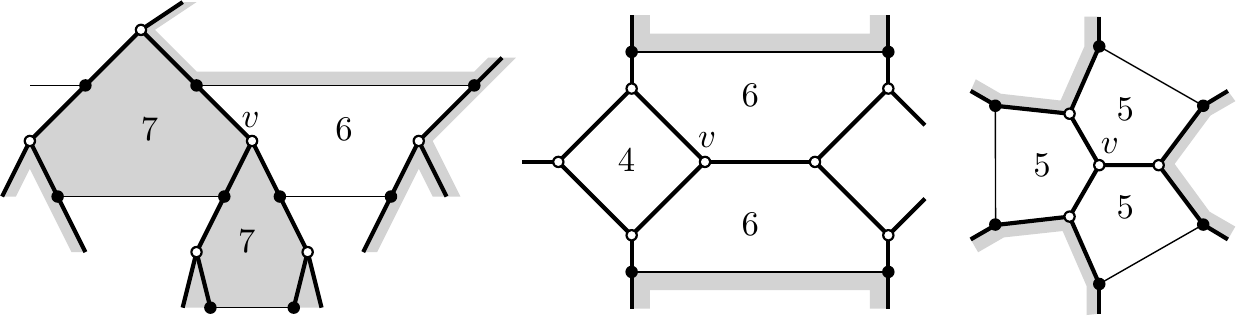}
\caption{From left to right: face sizes around $v$ as an inner vertex, leaf, and root, respectively.}
\label{fig:detail}
\end{figure}

The faces of $G_k$ have lengths $4,5,6$ or $7$. 
Indeed, every face contains a tree vertex (white) thus it is enough to analyze faces around all the tree vertices.
Let $v$ be a vertex of one of the underlying trees in $G_k$.
If $v$ is the root, then all faces adjacent to $v$ are of length $5$.
If $v$ is an internal vertex, then it is adjacent to a $6$-face
and two $7$-faces with an exception that when $v$ is adjacent to leaves of its tree, then one of the $7$-faces
degenerates to a $4$-face.
If $v$ is a leaf, then it is adjacent to a $4$-face and two $6$-faces.
See Figure~\ref{fig:detail}.

Clearly, $|V(G_k)|=n\geq 2^k$.
On the other hand, every two vertices in $G_k$ can be joined by a path of length at most $3k$.
Therefore the diameter of $G_k$ is at most $3k\leq 3\log (n)$.

\medskip

Let us finally mention that in~\cite{NS16} fullerene graphs of diameter $\sqrt{\frac{4n}{3}}$ have been constructed. Maybe this is the smallest diameter a fullerene graph can have.

\subsubsection*{Acknowledgments} The first author was supported by the  \emph{Agence nationale de la recherche} through project ANR-17-CE40-0015 and by the  \emph{Ministerio de Econom\'ia,
Industria y Competitividad} through grant RYC-2017-22701.

\bibliographystyle{my-siam}
{\footnotesize \bibliography{fullerrenes}}

\end{document}